\documentclass{article}
\usepackage[latin9]{inputenc}
\usepackage{amsmath}
\usepackage{amssymb}
\usepackage{esint}

\makeatletter

\usepackage{a4}\usepackage{amsthm}\usepackage{setspace}

\makeatother

\begin{document}

\title{Log-concavity of a Mixture of Beta Distributions%
\thanks{The author thanks Scott D. Kominers and Kareen Rozen for comments
and suggestions. %
}}

\author{Xiaosheng Mu\\
 Department of Economics\\
 Harvard University}

\maketitle
\doublespacing

\section{Introduction}

A non-negative function $f(x)$ defined on an interval $(a,b)$ is
said to be \emph{logarithmic concave} (\emph{log-concave}) if for
every $x,y\in(a,b)$ and every $0<\lambda<1$, we have 
\begin{equation}
f(\lambda x+(1-\lambda)y)\geq[f(x)]^{\lambda}[f(y)]^{1-\lambda}.
\end{equation}
If the inequality in (1) is reversed, the function $f$ is said to
be \emph{log-convex}. 

An equivalent definition of log-concavity (resp.~log-convexity) is
that the product $f(x)\cdot f(y)$ decreases (resp.~increases) in
$|x-y|$, holding $x+y$ fixed. Likewise, a non-negative sequence
$\{a_{i}\}_{i=0}^{n}$ is said to be log-concave if the product $a_{i}a_{j}$
decreases with $|i-j|$, holding $i+j$ fixed. It is known that log-concave
functions and sequences are closed under multiplication, integration,
and convolution \cite{prekopa}. \\

Proschan shows that \textit{log-convex} functions are closed under
addition or arbitrary mixture \cite{proschan}. However, the same
is not true for \textit{log-concave} functions, see examples in Barlow
and Proshcan \cite{barlow}. Except for the work of Lynch \cite{lynch}
and Block et al.~\cite{block}, little is known about general conditions
that guarantee the mixture of log-concave functions to also be log-concave. 

In this paper, we prove the following theorems: \newtheorem*{thm1}{Theorem
1} \begin{thm1} If $M>1$ and $\alpha(s)$ is a log-concave function
on $(0,M)$, then the function 
\begin{equation}
f(x)=\int_{0}^{M}\alpha(s)\binom{M}{s}(1-x)^{s}x^{M-s}\, ds
\end{equation}
is log-concave in $x$ on (0,1), where $\binom{M}{s}$ denotes $\frac{\Gamma(M+1)}{\Gamma(s+1)\Gamma(M-s+1)}=\frac{1}{(M+1)\beta(s+1,M-s+1)}$.
Here $\Gamma(\cdot)$ denotes the Gamma function, and $\beta(\cdot,\cdot)$
denotes the Beta function.\\

\end{thm1} \newtheorem*{thm2}{Theorem 2} \begin{thm2} If M
is a positive integer and $\{\alpha_{i}\}_{i=0}^{M}$ is a log-concave
sequence, then the function 
\begin{equation}
g(x)=\sum_{i=0}^{M}\alpha_{i}\binom{M}{i}(1-x)^{i}x^{M-i}
\end{equation}
is log-concave in x on the interval $(0,1)$. \\

\end{thm2} A direct corollary of these theorems is the following:
\newtheorem*{cor}{Corollary} \begin{cor} If $M>1$, then any
log-concave mixture of distributions 
\[
\{Beta(M-s+1,s+1)\}_{0<s<M}
\]
 has a log-concave density. Furthermore if $M$ is an integer, then
any discrete log-concave mixture of distributions 
\[
\{Beta(M-i+1,i+1)\}_{i=0}^{M}
\]
 also has a log-concave density. \end{cor}

While similar to the conditions given by Lynch in \cite{lynch}, this
result is not implied by those of Lynch because Beta densities are
not \textit{jointly} log-concave in the argument $x$ and the parameter
$s$. \\

\section{Proof of the Theorems}

Before proceeding to the proofs, we recall two technical lemmata that
will be of use: \newtheorem*{lemma1}{Lemma 1} \begin{lemma1}
Let $a(q),b(q),u(q),v(q)$ be non-negative continuous functions defined
on the interval $[0,m]$, such that $a(q)$ is decreasing in $q$.
Suppose further that for each $q$, $a(q)\geq b(q)\geq0$, and $\int_{r=0}^{q}u(r)\, dr\geq\int_{r=0}^{q}v(r)\, dr$.
Then we have 
\begin{align*}
\int_{q=0}^{m}a(q)u(q)\, dq\geq\int_{q=0}^{m}b(q)v(q)\, dq.
\end{align*}
\end{lemma1} \proof Define $U(q)=\int_{r=0}^{q}u(r)\, dt$ and $V(q)=\int_{r=0}^{q}v(r)\, dr$.
Using integration by parts, we have 
\[
\int_{q=0}^{m}a(q)u(q)\, ds=\int_{q=0}^{m}a(q)\, dU(q)=a(m)U(m)+\int_{q=0}^{m}U(q)\, d(-a(q)).
\]

Since $-a(q)$ is increasing and $U(q)\geq V(q)$ pointwise, the Stieltjes
integral $\int_{q=0}^{m}U(q)\, d(-a(q))$ is greater than or equal
to $\int_{q=0}^{m}V(q)\, d(-a(q))$. Another use of integration by
parts shows that 
\[
\int_{q=0}^{m}a(q)u(q)\, dq\geq a(m)V(m)+\int_{q=0}^{m}V(q)\, d(-a(q))=\int_{q=0}^{m}a(q)v(q)\, dq\geq\int_{q=0}^{m}b(q)v(q)\, dq.
\]

In the literature, Lemma 1 is often called the ``majorization trick''
(see for instance \cite{maj}). \\

\newtheorem*{lemma2}{Lemma 2} \begin{lemma2} Let $M>1,q>0$
and $n>-2$ be fixed parameters. Define the sets: 

\[
\begin{array}{c}
A=\{s:|s-(n-s)|\leq q,0\leq s,n-s\leq M\};\\
B=\{s:|s-(n+1-s)|\leq q,0\leq s,n+1-s\leq M\};\\
C=\{s:|s-(n-s)|\leq q,0\leq s+1,n+1-s\leq M\}.
\end{array}
\]
Then the following inequalities hold: 
\begin{align}
\int_{A}\binom{M-1}{s}\binom{M-1}{n-s}\, ds & \geq\int_{A}\binom{M}{s}\binom{M-2}{n-s}\, ds;\\
\int_{B}\binom{M-1}{s}\binom{M-1}{n-s}\, ds & \leq\int_{B}\binom{M}{s}\binom{M-2}{n-s}\, ds;\\
\int_{C}\binom{M-1}{s}\binom{M-1}{n-s}\, ds & \geq\int_{C}\binom{M}{s+1}\binom{M-2}{n-s-1}\, ds.
\end{align}
\end{lemma2} \proof We only give a proof for (4), as (5) and (6)
follow along similar lines. Using properties of the $\Gamma$ function---in
particular $\Gamma(M)=(M-1)\cdot\Gamma(M-1),$---we see that the generalized
binomial coefficients satisfy 
\begin{equation}
\begin{split}\binom{M}{s}=\binom{M}{M-s}>0 & \text{, }\forall-1<s<M+1;\\
\binom{M-1}{s}=\frac{M-s}{M}\binom{M}{s} & \text{, }\binom{M-1}{s-1}=\frac{s}{M}\binom{M}{s};\\
\binom{M}{s}=\binom{M-1}{s} & +\binom{M-1}{s-1}.
\end{split}
\end{equation}

To prove (4), first note that we can assume $n>0$ and $q\leq\min{\{n,2M-n\}}$,
which is the maximum difference between $s$ and $n-s$ when $s\in A$.
Writing $k=\frac{n-q}{2}$, we have $A=[k,n-k]\subset[0,M]$. Using
(7) to write $\binom{M-1}{n-s}=\binom{M-2}{n-s}+\binom{M-2}{n-s-1}$
and $\binom{M}{s}=\binom{M-1}{s}+\binom{M-1}{s-1}$, we can calculate
the difference between the two sides of (4) as: 
\begin{equation}
\begin{split} & \int_{s=k}^{n-k}\left(\binom{M-1}{s}\binom{M-2}{n-s-1}-\binom{M-1}{s-1}\binom{M-2}{n-s}\right)\, ds\\
= & \int_{s=k-1}^{k}\left(\binom{M-1}{n-s-1}\binom{M-2}{s}-\binom{M-1}{s}\binom{M-2}{n-s-1}\right)\, ds\\
= & \int_{s=k-1}^{k}\frac{n-2s-1}{M-1}\binom{M-1}{s}\binom{M-1}{n-s-1}\, ds.
\end{split}
\end{equation}

By (7), the two binomial coefficients above are non-negative in the
range of integration. When $k\leq\frac{n-1}{2}$, the term $n-2s-1$
is always non-negative, so is the last line of (8). When $k>\frac{n-1}{2}$,
we can write the last line of (8) as 
\begin{equation}
\int_{s=k-1}^{n-k-1}\frac{n-2s-1}{M-1}\binom{M-1}{s}\binom{M-1}{n-s-1}\, ds+\int_{s=n-k-1}^{k}\frac{n-2s-1}{M-1}\binom{M-1}{s}\binom{M-1}{n-s-1}\, ds.
\end{equation}

Note that the integrand $\frac{n-2s-1}{M-1}\binom{M-1}{s}\binom{M-1}{n-s-1}$
is non-negative when $s\leq n-k-1\leq\frac{{n-1}}{2}$, and it is
an odd function with respect to $s=\frac{{n-1}}{2}.$ Thus the first
integral in (9) is non-negative, while the second evaluates to zero.
(8) follows, so does the lemma. \\

For the sake of completeness we state below the discrete analogs of
the preceding two lemmata, which will be used in the proof of Theorem
2. We omit the proofs because they are the same. \newtheorem*{lemma1'}{Lemma
1'} \begin{lemma1'} Let $a_{1}\geq a_{2}\geq\dots\geq a_{m}\geq0,b_{1},b_{2},\dots,b_{m}\geq0$
be two sequences of real numbers satisfying $a_{q}\geq b_{q}$ for
each $q$. Consider two more sequences of non-negative real numbers
$u_{1},u_{2},\dots u_{n}$ and $v_{1},v_{2},\dots v_{n}$, such that
$\sum_{r=1}^{q}u_{r}\geq\sum_{r=1}^{q}v_{r}$ for each $q$. Then,
\[
\sum_{q=1}^{m}a_{q}u_{q}\geq\sum_{q=1}^{m}b_{q}v_{q}.
\]
\end{lemma1'}

\newtheorem*{lemma2'}{Lemma 2'} \begin{lemma2'} If $M$ is a
positive integer, $n$ is a non-negative integer and $k$ is an integer
such that $k\leq\frac{n+1}{2}$, then the following inequalities hold:
\begin{align*}
\sum_{i=k}^{n-k}\binom{M-1}{i}\binom{M-1}{n-i} & \geq\sum_{i=k}^{n-k}\binom{M}{i}\binom{M-2}{n-i};\\
\sum_{i=k}^{n-k+1}\binom{M-1}{i}\binom{M-1}{n-i} & \leq\sum_{i=k}^{n-k+1}\binom{M}{i}\binom{M-2}{n-i};\\
\sum_{i=k}^{n-k}\binom{M-1}{i}\binom{M-1}{n-i} & \geq\sum_{i=k}^{n-k}\binom{M}{i+1}\binom{M-2}{n-i-1};
\end{align*}
where as usual we define $\binom{M}{i}=0$ when $i<0$ or $i>M$.
\\
 \end{lemma2'}

\noindent \textbf{Proof of Theorem 1:} We first make some preliminary
simplifications. When $\alpha(s)$ is identically zero, the result
is trivial. Otherwise 
\[
f(x)=\int_{s=0}^{M}\alpha(s)\binom{M}{s}(1-x)^{s}x^{M-s}\, ds
\]
 is strictly positive for $x\in(0,1)$. Thus $\log{f}$ is well-defined
on the open interval. Its derivative is $\frac{f'}{f}$ and its second
derivative is $\frac{f'^{2}-f\cdot f''}{f^{2}}$. It thus suffices
to show that $f'(x)^{2}\geq f(x)\cdot f''(x)$. 

We will show that the following stronger inequality holds:

\begin{equation}
\frac{M-1}{M}f'(x)^{2}\geq f(x)\cdot f''(x)\text{ for every \ensuremath{x} in (0,1). }
\end{equation}

Using (7), the derivative of $\binom{M}{s}(1-x)^{s}x^{M-s}$ is 
\[
\begin{split} & \binom{M}{s}(-s)(1-x)^{s-1}x^{M-s}+\binom{M}{s}(1-x)^{s}(M-i)x^{M-s-1}\\
= & M\left(-\binom{M-1}{s-1}(1-x)^{s-1}x^{M-s}+\binom{M-1}{s}(1-x)^{s}x^{M-s-1}\right).
\end{split}
\]

It follows that 
\[
\begin{split} & f'(x)=M\left(\int_{s=-1}^{M}[\alpha(s)-\alpha(s+1)]\binom{M-1}{s}(1-x)^{s}x^{M-s-1}\, ds\right);\\
 & f''(x)=M(M-1)\left(\int_{s=-2}^{M}[\alpha(s)-2\alpha(s+1)+\alpha(s+2)]\binom{M-2}{s}(1-x)^{s}x^{M-s-2}\, ds\right);
\end{split}
\]
where we define $\alpha(s)=0$ whenever $s<0$ or $s>M$.\\

Introducing another dummy variable $t$, we turn the desired inequality
(10) into: 
\begin{equation}
\begin{split} & \left(\int_{-1}^{M}[\alpha(s)-\alpha(s+1)]\binom{M-1}{s}(1-x)^{s}x^{M-s-1}\, ds\right)\left(\int_{-1}^{M}[\alpha(t)-\alpha(t+1)]\binom{M-1}{t}(1-x)^{t}x^{M-t-1}\, dt\right)\\
 & \geq\left(\int_{0}^{M}\alpha(s)\binom{M}{s}(1-x)^{s}x^{M-s}\, ds\right)\left(\int_{-2}^{M}[\alpha(t)-2\alpha(t+1)+\alpha(t+2)]\binom{M-2}{t}(1-x)^{t}x^{M-t-2}\, dt\right)
\end{split}
\end{equation}

By expanding the products on both sides of (11) and collecting terms
that have $(1-x)^{n}x^{2M-n-2}$ in common, it suffices to prove that
for any $n$ $(-2<n<2M-2)$, 
\begin{equation}
\begin{split} & \int_{s+t=n}[\alpha(s)-\alpha(s+1)][\alpha(t)-\alpha(t+1)]\binom{M-1}{s}\binom{M-1}{t}\, ds\\
\geq & \int_{s+t=n}\alpha(s)[\alpha(t)-2\alpha(t+1)+\alpha(t+2)]\binom{M}{s}\binom{M-2}{t}\, ds.
\end{split}
\end{equation}

We show that the following three inequalities hold, which collectively
imply (12): 
\begin{equation}
\int_{s+t=n}\alpha(s)\alpha(t)\binom{M-1}{s}\binom{M-1}{t}\, ds\geq\int_{s+t=n}\alpha(s)\alpha(t)\binom{M}{s}\binom{M-2}{t}\, ds;
\end{equation}
\begin{equation}
\int_{s+t=n}\alpha(s)\alpha(t+1)\binom{M-1}{s}\binom{M-1}{t}\, ds\leq\int_{s+t=n}\alpha(s)\alpha(t+1)\binom{M}{s}\binom{M-2}{t}\, ds;
\end{equation}
\begin{equation}
\int_{s+t=n}\alpha(s+1)\alpha(t+1)\binom{M-1}{s}\binom{M-1}{t}\, ds\geq\int_{s+t=n}\alpha(s)\alpha(t+2)\binom{M}{s}\binom{M-2}{t}\, ds.
\end{equation}

To prove (13), we use the ``majorization trick'' of Lemma 1. 

Let $a(q)=b(q)=\alpha(\frac{n-q}{2})\cdot\alpha(\frac{n+q}{2})$,
then the log-concavity of $\alpha(\cdot)$ ensures that $a(q)$ is
non-negative and decreasing in $q$. Moreover let 
\[
u(q)=2\binom{M-1}{\frac{n-q}{2}}\binom{M-1}{\frac{n+q}{2}},\text{ }v(q)=\binom{M}{\frac{n-q}{2}}\binom{M-2}{\frac{n+q}{2}}+\binom{M-2}{\frac{n-q}{2}}\binom{M}{\frac{n+q}{2}}.
\]

Inequality (13) becomes $\int_{q=0}^{n}a(q)u(q)\, dq\geq\int_{q=0}^{n}b(q)v(q)\, dq.$
It remains to check the majorization condition 
\[
\int_{r=0}^{q}u(r)\, dr\geq\int_{r=0}^{q}v(r)\, dr\text{ for each q.}
\]
But this is exactly the inequality (4) in Lemma 2. So (13) is proved.
In the same way, (14) and (15) reduce to (5) and (6) in Lemma 2. Theorem
1 follows. \\

We note that by (7), the function $\binom{M}{s}(1-x)^{s}x^{M-s}$
is well defined and positive even for $-1<s<0$ and $M<s<M+1.$ However,
it fails to be log-concave on these intervals. Thus Theorem 1 would
fail if we allow the range of integration in (2) to include these
intervals. \\

\noindent \textbf{Proof of Theorem 2:} We take the same steps as before
to reduce the result to inequalities (13), (14) and (15), with the
integral replaced by a discrete sum. Those inequalities are in turn
implied by Lemma 1' and Lemma 2'. The theorem follows.\\

It is worth noting that the discrete analog of inequality (10) is
tight. In fact, whenever the $\alpha$ sequence is geometric, we have
\[
g(x)=\sum_{i=0}^{M}\alpha_{i}\binom{M}{i}(1-x)^{i}x^{M-i}=c(1+\lambda x)^{M}
\]
for some $c>0$ and $\lambda$. Thus $g(x)\cdot g''(x)=\frac{M-1}{M}g'(x)^{2}$,
which is slightly less than $g'(x)^{2}$ for large $M$. This suggests
that improvement of Theorem 2 (likewise Theorem 1) would depend on
a different type of conditions on the mixing sequence/function $\{\alpha_{i}\}_{i=0}^{M}$,
rather than its own log-concavity. \\

\noindent \textbf{Proof of the Corollary: } Recall that the $Beta(M-s+1,s+1)$
distribution has density 
\[
\frac{1}{\beta(M-s+1,s+1)}x^{M-s}(1-x)^{s}=(M+1)\binom{M}{s}(1-x)^{s}x^{M-s}.
\]

Therefore a mixture Beta distribution has density given by the form
in Theorems 1 and 2. From the theorems, this density is log-concave
whenever the weights are log-concave.

\section{Applications}

A recent application of this paper is seen in \cite{rozen}, where
the authors assume subpopulations of consumers having distinct attention
capacity. With a total of $M$ markets, the probability that a particular
consumer pays attention to a market takes the form of a Beta distribution
function in the price charged. The expected amount of attention that
the market receives is therefore a mixture of Beta distributions.
Using Theorem 2, the authors characterize equilibrium strategies and
establish that they are monotone \cite{rozen}. Log-concavity is important
elsewhere in economics, see Bagnoli and Bergstrom \cite{econ}.\\

In the statistical theory of reliability, one is interested in the
failure rate of mixture systems. The techniques here thus complement
earlier work by Lynch \cite{lynch} and Block et al.~\cite{block}.
Finally, log-concavity has powerful implications for truncated distributions
\cite{trun}, hypothesis testing (e.g. the Karlin-Rubin theorem \cite{test}),
and maximum likelihood estimation \cite{mle}. The recent surge of
interest in mixture models renders it necessary to gain further understanding
of the log-concavity of general mixture distributions \cite{spec},
\cite{clus}. This paper is a step toward that understanding for the
special family of Beta distributions. \\

\end{document}